\documentclass[letterpaper,12pt]{article}
\usepackage{amssymb}
\usepackage{amsmath}
\usepackage{makeidx}
\usepackage{amscd}
\usepackage[arrow,matrix]{xy}
\usepackage{hyperref}

\newtheorem{theo}{Theorem}[section]
\newtheorem{prop}[theo]{Proposition}
\newtheorem{lem}[theo]{Lemma}
\newtheorem{cor}[theo]{Corollary}
\newtheorem{defi}[theo]{Definition}

\def \q {{\mathfrak q}}

\def \K{{\mathcal K}}

\def \Br {{\rm{Br}}}

\def \Ga {{\Gamma}}
\def \R {{\mathbb{R}}}
\def \Pic {{\rm {Pic}}}

\def \Gal {{\rm{Gal}}}
\def \Ker {{\rm{Ker}}}

\def \Frob{{\rm{Frob}}}

\def \A{{\mathbb A}}
\def \P{{\mathbb P}}
\def \Spec {{\rm{Spec}}}
\def \dim {{\rm{dim}}}
\def \Hom {{\rm {Hom}}}
\def \End {{\rm {End}}}
\def \Pic {{\rm {Pic}}}
\def \GL {{\rm {GL}}}
\def \SL {{\rm {SL}}}

\def \Aut{{\rm Aut}}

\def\ov{\overline}

\def \Z {{\mathbb Z}}
\def \Q {{\mathbb Q}}
\def \F {{\mathbb F}}

\def \GG {{\rm G}}

\def \Tr {{\rm{Tr}}}

\def \val {{\rm{val}}}

\def\G{{\mathbb G}}

\def\GG{{\cal G}}

\def\sC{{\mathcal C}}
\def\sA{{\mathcal A}}

\def\lra{\longrightarrow}

\def\N{{\rm N}}
\def\H{{\rm H}}

\def\Tr{{\rm Tr}}
\def\Aff{{\rm Aff}}
\def\O{{\cal O}}

\def\inv{{\rm inv}}
\def\Kum{{\rm Kum}}

\def\NS{{\rm NS\,}}

\def\O{{\cal O}}

\def\val{{\rm val}}
\def\Sel{{\rm Sel}}

\def\discr{{\rm discr}}
\def\Ga{\Gamma}

\def\loc{{\rm loc}}

\def\m{{\mathfrak m}}

\newcommand{\bthe}{\begin{theo}}
\newcommand{\ble}{\begin{lem}}
\newcommand{\bpr}{\begin{prop}}
\newcommand{\bco}{\begin{cor}}
\newcommand{\bde}{\begin{defi}}
\newcommand{\ethe}{\end{theo}}
\newcommand{\ele}{\end{lem}}
\newcommand{\epr}{\end{prop}}
\newcommand{\eco}{\end{cor}}
\newcommand{\ede}{\end{defi}}

\DeclareFontEncoding{OT2}{}{} 
\DeclareTextFontCommand{\textcyr}{\fontencoding{OT2}\fontfamily{wncyr}\fontseries{m}\fontshape{n}\selectfont}

\newcommand{\Sha}{\textcyr{Sh}}

\title{Hasse principle for generalised Kummer varieties}
\author{Yonatan Harpaz and Alexei N. Skorobogatov}
\date{\today}
\begin{document}
\baselineskip=15pt
\maketitle

\begin{abstract}
\noindent 
The existence of rational points on Kummer varieties associated to $2$-coverings
of abelian varieties over number fields
can sometimes be proved through the variation of the Selmer group in
the family of quadratic twists of the underlying abelian variety,
using an idea of Swinnerton-Dyer. Following Mazur and Rubin,
we consider the case when the Galois
action on the 2-torsion has a large image. Under mild additional hypotheses
we prove the Hasse principle for the associated Kummer varieties 
assuming the finiteness of relevant Shafarevich-Tate groups.
\end{abstract}

\section{Introduction}

The aim of this paper is to give some evidence in favour of the conjecture
that the Brauer--Manin obstruction is the only obstruction to the Hasse
principle for rational points on K3 surfaces over number fields, see
\cite[p. 77]{SkOb} and \cite[p. 484]{SZ}. Conditionally on the 
finiteness of relevant Shafarevich--Tate groups we establish the Hasse
principle for certain families of Kummer surfaces. These 
surfaces are quotients of 2-coverings of an abelian surface $A$ 
by the antipodal involution, where

(1) $A$ is the product of elliptic curves $A=E_1\times E_2$, or

(2) $A$ is the Jacobian of a curve $C$ of genus 2 with a 
rational Weierstra{\ss} point.

\noindent Both cases are treated by the same method which allows us
to prove a more general result for generalised Kummer varieties $\Kum(Y)$,
where $Y$ is a 2-covering of an abelian variety $A$ over a number field $k$,
provided certain conditions are satisfied. By a 2-covering
we understand a torsor $Y$ for $A$ 
such that the class $[Y]\in\H^1(k,A)$ has order at most 2, so that 
$Y$ is the twist of $A$ by a 1-cocycle with
coefficients in $A[2]$. Then the antipodal
involution $\iota:A\to A$ acts on $Y$ and we can define $\Kum(Y)$
as the minimal desingularisation of $Y/\iota$. In this
introduction we explain the results pertaining to cases (1) and (2)
above and postpone the statement of a more general theorem 
until the next section.
In case (1) we have the following result whose proof can be found at 
the end of Section \ref{2}.

\medskip

\noindent{\bf Theorem A}
{\em Let $g_1(x)$ and $g_2(x)$ be irreducible polynomials of degree $4$ 
over a number field $k$, each with the Galois group $S_4$. 
Let $w_1$ and $w_2$ be distinct odd places of $k$
such that for all $i,j\in \{1,2\}$
the coefficients of $g_i(x)$ are integral at $w_j$
and $\val_{w_j}(\discr(g_i))=\delta_{ij}$. Assume the finiteness
of the $2$-primary torsion subgroup of the Shafarevich--Tate group for 
each quadratic twist
of the Jacobian of the curve $y^2=g_i(x)$, where $i=1,2$. If the
Kummer surface with the affine equation $z^2=g_1(x)g_2(y)$
is everywhere locally soluble, then it has a Zariski dense
set of $k$-points.}

\medskip

The conditions of Theorem A probably hold for `most' pairs of
polynomials of degree 4 with coefficients in the
ring of integers $\O_k$ of $k$. Indeed,
by a theorem of van der Waerden (see \cite[Thm. 1]{Coh} for a statement
over an arbitrary number field $k$)
100 \% of such polynomials have the Galois group $S_4$,
when ordered by the height of coefficients.
By a recent theorem of M. Bhargava, C. Skinner and W. Zhang \cite[Thm. 2]{BSZ}
a majority of elliptic curves over $\Q$,
when ordered by na\"ive height, have finite Shafarevich--Tate groups.
An easy sieve argument shown to us by Andrew Granville
gives that for fixed rational primes $w_1$ and $w_2$
the condition on the discriminants of $g_1$ and $g_2$
is satisfied for a positive proportion of pairs of 
polynomials of degree 4 in $\Z[x]$.

To give an explicit description of our results in case (2) 
we need to recall the realisation of Kummer surfaces attached to the 
Jacobian of a genus 2 curve as complete intersections
of three quadrics in $\P^5_k$. We mostly follow \cite{Sk10}; 
for the classical theory over an algebraically
closed field see \cite[Ch. 10]{D_CAG}.

Let $C$ be the hyperelliptic curve $y^2=f(x)$, 
where $f(x)$ is a monic separable polynomial
of degree 5 over a field $k$ of characteristic different from 2.
Let $L$ be the \'etale $k$-algebra $k[x]/(f(x))$, and let
$\theta\in L$ be the image of $x$. Let $A$ be the Jacobian of $C$.
Then $A[2]={\rm R}_{L/k}(\mu_2)/\mu_2$, where ${\rm R}_{L/k}$ is the Weil
restriction of scalars. It follows that $\H^1(k,A[2])=L^*/k^*L^{*2}$.
The $k$-torsor for $A[2]$ defined by $\lambda\in L^*$ is the closed subset
of ${\rm R}_{L/k}(\G_m)/\{\pm 1\}$ given by $\lambda=z^2$.
Let $Y_\lambda$ be the 2-covering of $A$
obtained by twisting $A$ by the 1-cocycle
represented by $\lambda\in L^*$. Then $\Kum(Y_\lambda)$ is a
subvariety of $\P({\rm R}_{L/k}(\A^1_L)\times \A^1_k)\simeq\P^5_k$
given by three quadratic equations
\begin{equation}\Tr_{L/k}\left(\lambda \frac{u^2}{f'(\theta)}\right)=
\Tr_{L/k}\left(\lambda \frac{\theta u^2}{f'(\theta)}\right)=
\Tr_{L/k}\left(\lambda \frac{\theta^2u^2}{f'(\theta)}\right)-{\rm N}_{L/k}(\lambda)u_0^2=0,
\label{kummer}\end{equation}
where $u$ is an $L$-variable, $u_0$ is a $k$-variable, and $f'(x)$ 
is the derivative of $f(x)$. If $\lambda\in k^*L^{*2}$,
then an easy change of variable reduces (\ref{kummer}) to 
the same system of equations with $\lambda=1$. The Lagrange interpolation
identities imply that $\Kum(Y_1)$ contains the projective line
given parametrically by $u=r+s\theta$, $u_0=s$, so for the purpose of 
establishing the Hasse principle this case can be excluded.

\medskip

\noindent{\bf Theorem B}
{\em Let $f(x)$ be a monic irreducible polynomial of degree $5$
over a number field $k$, and let $L=k[x]/(f(x))$.
Let $w$ be an odd place of $k$ such that the coefficients of $f(x)$ 
are integral at $w$ and $\val_w(\discr(f))=1$. 
Assume the finiteness of the $2$-primary torsion subgroup of the
Shafarevich--Tate group for each quadratic twist
of the Jacobian of the curve $y^2=f(x)$. 
Let $\lambda\in L^*$ be such that for some $r\in k^*$ the valuation
of $\lambda r$ at each completion of $L$ over $w$ is even, but
$\lambda\notin k^*L^{*2}$. If the Kummer surface given by $(\ref{kummer})$
is everywhere locally soluble,
then it has a Zariski dense set of $k$-points.}

\medskip

See the end of the next section for the proof.
Any Kummer surface (\ref{kummer}) can be mapped to $\P_k^3$
by a birational morphism that contracts 16 disjoint rational curves
onto singular points. This map is defined by the linear system 
on $Y_\lambda$ obtained from the linear system $|2C|$
on $A$. The image of $Y_\lambda$ is a quartic surface
$X\subset\P^3_k$ which is the classical singular Kummer surface
with 16 nodes. (See \cite{D_CAG} and \cite{GD} for a modern
account of the geometry of $X$ over an algebraically closed field.)
The group of projective automorphisms
of $X$ is $A[2]$ and the singular locus $X_{\rm sing}$ is 
a $k$-torsor for $A[2]$ such that $X$ is obtained by
twisting $A/\iota$ by $X_{\rm sing}$. The condition
$\lambda\notin k^*L^{*2}$, which we need to prove the Zariski density
of $X(k)$, is precisely the condition that this torsor is non-trivial,
or, equivalently, $X_{\rm sing}\cap X(k)=\emptyset$.

The main idea of the proof of Theorems A and B is due to Swinnertion-Dyer.
Let $\alpha\in\H^1(k,A[2])$ be the class of a 1-cocycle used to obtain
$Y$ from $A$. For an extension $F/k$ of degree at most 2,
we denote by $A^F$ and $Y^F$ the corresponding quadratic twists of $A$ 
and $Y$, respectively. Then $Y^F$ is a torsor for $A^F$
defined by the same $\alpha\in \H^1(k,A^F[2])=\H^1(k,A[2])$.
To find a rational point on $\Kum(Y)$ it is enough to find a rational
point on $Y^F$ for some $F$. At the first step of the proof,
using a fibration argument one produces
a quadratic extension $F$ such that $Y^F$ is everywhere locally soluble.
Equivalently, $\alpha\in \H^1(k,A^F)$ is in the 2-Selmer group of $A^F$.
At the second step one constructs $F$ such that the 2-Selmer group of $A^F$
is spanned by $\alpha$ and the image of $A^F[2](k)$ under the Kummer 
map. (In the cases considered in this paper $A^F[2](k)=A[2](k)=0$.)
This implies that $\Sha(A^F)[2]$ is $\Z/2$ or $0$. 
In previous applications of the method \cite{SD2001, SSD}
$A$ was a product of two elliptic curves, in which case 
the Cassels--Tate pairing on $\Sha(A^F)$ is alternating. 
The assumption that $\Sha(A^F)$ is finite then implies that 
the order of $\Sha(A^F)[2]$ is a square and hence 
$\Sha(A^F)[2]=0$. In particular, $Y^F\simeq A^F$ has a $k$-point.
In this paper we consider more general principally polarised abelian varieties.
The theory developed by Poonen and 
Stoll in~\cite{PS} ensures that in the cases considered here
the Cassels--Tate pairing on $\Sha(A^F)$ defined using
the principal polarisation is still alternating, so the proof can be
concluded as before.

Swinnertion-Dyer's method was used in combination with Schinzel's 
Hypothesis (H) in \cite{CSS, SD2000, W}. The immediate precursor 
of our Theorem A is \cite{SSD}, which treats Kummer surfaces attached 
to products of elliptic curves without assuming Schinzel's Hypothesis
(see also \cite{SD2001} for a case with complex multiplication). 
Central to Swinnertion-Dyer's method is a linear algebra
construction that represents the Selmer group as the kernel of a 
symmetric bilinear form. The difficulty of operating this machinery
makes implementation of the method a rather delicate task requiring a
considerable effort.
In the present paper this linear algebra machinery is not used.
Instead we use the ideas of Mazur and Rubin from \cite{MR07} and especially
from \cite{MR10}. 

Let us note that given an elliptic curve
$E$ over a number field $k$ it is not always possible to find 
a quadratic extension $F/k$ such that
the 2-Selmer group of $E^F$ is spanned by a fixed class 
$\alpha\in\H^1(k,E[2])$ and the image of $E^F[2](k)$. 
Firstly, the parity of
the rank of the 2-Selmer group of $E^F$ can be the same for all $F$:
this happens precisely when $k$ is totally imaginary and $E$
acquires everywhere good reduction over an abelian extension
of $k$, see \cite[Remark 4.9]{DD}. Secondly, 
over any number field $k$ there are
elliptic curves $E$ such that for any quadratic extension
$F/k$ the difference between
the 2-Selmer rank of $E^F$ and the dimension of the $\F_2$-vector space
$E[2](k)$ is at least the number of complex places of $k$,
see \cite{ZK1, ZK2}. Such examples can occur when $E[2](k)\simeq\Z/2$
and $E$ has a cyclic isogeny of degree 4 defined over $k(E[2])$
but not over $k$.

In this paper we do not
discuss the conjecture that rational points of a K3 surface
are dense in its Brauer--Manin set \cite[p. 77]{SkOb}, \cite[p. 484]{SZ}.
A recent result of D. Holmes and R. Pannekoek \cite{HP} 
shows that if this 
conjecture is extended to all generalised Kummer varieties, then
the ranks of quadratic twists of any given abelian variety
over a given number field are not bounded.

The main technical result of the paper is Theorem \ref{t:main}.
It is stated in Section \ref{2} where we also show that 
Theorem \ref{t:main} implies Theorems A and B. 
The Brauer--Manin obstruction does not appear in these results
for purely algebraic reasons, see the remarks after Theorem \ref{t:main}.
In Section \ref{one} we systematically develop
the Galois-theoretic aspect of the approach of Mazur and Rubin.
We recall the necessary facts about the Kummer
map for quadratic twists of abelian varieties 
over local fields in Section \ref{local}. In Section \ref{global}
we discuss the Selmer group and the Cassels--Tate pairing over
a number field. A reduction to everywhere soluble 2-coverings
is carried out in Section \ref{fibr} using a known case of the
fibration method. We finish the proof of Theorem \ref{t:main} in Section
\ref{proof}.

The second named author is grateful to the Radboud
University in Nijmegen for its hospitality. We would like
to thank Tim Dokchitser, 
Evis Ieronymou and Ren\'e Pannekoek for helpful discussions.

\section{Main results} \label{2}

Let $k$ be a field of characteristic different from 2
with a separable closure $\bar k$ and the Galois group
$\Ga_k=\Gal(\bar k/k)$.

Let $A$ be an abelian variety over $k$.
Let $K=k(A[2])\subset\bar k$ be the field of definition of $A[2]$, 
that is, the smallest field such that $A[2](K)=A[2](\bar k)$.
Let $G=\Gal(K/k)$.
Consider the following conditions:

\smallskip

{\bf (a)} $A[2]$ is a simple $G$-module and $\End_G(A[2])=\F_2$;

{\bf (b)} $\H^1(G,A[2])=0$;

{\bf (c)} there exists $g\in G$ such that $A[2]/(g-1)=\F_2$;

{\bf (d)} there exists $h\in G$ such that $A[2]/(h-1)=0$.

\ble \label{k1}
Let $A$ be the Jacobian of a smooth projective
curve with the affine equation $y^2=f(x)$,
where $f(x)\in k[x]$ is an irreducible polynomial of odd degree $m\geq 3$.
If the Galois group of $f(x)$ is the symmetric group on $m$ letters $S_m$,
then $A$ satisfies conditions (a), (b), (c), (d).
\ele
{\em Proof.} It is well known that the $\Ga_k$-module $A[2]$
is the zero-sum submodule of the vector space
$(\F_2)^m$ freely generated by the roots of $f(x)=0$
with the natural permutation action of $\Ga_k$. Since
$m$ is odd, the permutation $\Ga_k$-module $(\F_2)^m$ is the direct sum of
$A[2]$ and the $\F_2$-vector space spanned by the vector $(1,\ldots,1)$.

If an $S_m$-submodule of $(\F_2)^m$ contains a vector with at least one 
coordinate $0$ and at least one coordinate $1$, then it contains the zero-sum
submodule. Hence $A[2]$ is a simple $S_m$-module.
A direct calculation with matrices shows that the $m\times m$ matrices
commuting with all permutation matrices are the 
linear combinations of the identity and the all-1 matrix. 
We deduce that $\End_{S_m}(A[2])=\F_2$, thus (a) holds. 

The permutation $S_m$-module $(\F_2)^m$ is isomorphic to $\F_2[S_m/S_{m-1}]$.
By Shapiro's lemma we have 
$$\H^1(S_m, \F_2[S_m/S_{m-1}])=\H^1(S_{m-1},\F_2)=\Hom(S_{m-1},\F_2)=\F_2.$$
Since $\H^1(S_m,\F_2)=\F_2$, we obtain $\H^1(S_m,A[2])=0$, so (b) holds.

If $g$ is a cycle of length $m-1$, then $A[2]/(g-1)=\F_2$, so (c) is satisfied.
If $h$ is a cycle of length $m$, then $A[2]/(h-1)=0$, so 
(d) is satisfied. $\Box$

\medskip

Let $A_1,\ldots, A_r$ be abelian varieties over $k$.
For each $i=1,\ldots,r$ let $K_i=k(A_i[2])$ and $G_i=\Gal(K_i/k)$. 
We assume the following condition.

\smallskip

{\bf (e)} The fields $K_1,\ldots,K_r$ are linearly disjoint over $k$.

\smallskip

\noindent By definition this means that 
$[K_1\ldots K_r:k]=\prod_{i=1}^r [K_i:k]$. Thus the Galois group 
of $K_1\ldots K_r$ over $k$ is $\prod_{i=1}^r G_i$.

When $k$ is a {\em number field} we shall also assume the
following condition.

\smallskip

{\bf (f)} There exist distinct odd primes $w_1,\ldots, w_r$ of $k$
such that for each $i=1,\ldots, r$ the abelian variety
$A_i$ has bad reduction at $w_i$ and the number of 
geometric connected components of the N\'eron model
of $A_i$ at $w_i$ is odd, whereas each $A_j$ for $j\not=i$
has good reduction at $w_i$. 

\smallskip

Let $k_i^{\rm ab}$ be the maximal abelian subextension of $k\subset K_i$.
Equivalently, 
$\Gal(k_i^{\rm ab}/k)$ is the maximal abelian quotient $G_i^{\rm ab}$ of
$G_i$. Let us finally assume the condition

\smallskip

{\bf (g)} for each $i=1,\ldots, r$ the field $k_i^{\rm ab}$ does not contain
a non-trivial subextension of $k$ unramified at $w_i$. (In other words,
if $k\subset k'\subset k_i^{\rm ab}$ is a subfield such that $k'/k$
is unramified at $w_i$, then $k'=k$.) Equivalently, $G_i^{\rm ab}$ coincides
with the inertia subgroup of the unique prime of $k_i^{\rm ab}$ 
dividing $w_i$.

\smallskip

Let $F$ be a field extension of $k$ of degree at most 2. 
Associating to $F$ its quadratic character defines a bijection between 
the set of such extensions and $\H^1(k,\Z/2)=k^*/k^{*2}$.
Any abelian variety is endowed with an action of $\Z/2$
such that the generator acts as the antipodal involution $\iota(x)=-x$.
We denote by $A^F$ the quadratic twist of $A$ by $F$, that is,
the abelian variety over $k$ obtained
by twisting $A$ by the quadratic character of $F$ with respect to this
action of $\Z/2$.
For example, if $A$ is an elliptic curve with
the Weierstra{\ss} equation $y^2=f(x)$, then $A^F$ is
given by $y^2=cf(x)$, where $c\in k^*$ is such that $F=k(\sqrt{c})$.

Now we are ready to state the main theorem of this paper.
Recall that a class in $\H^1(k,A[2])$ is {\em unramified} at a
non-Archimedean place $v$ of $k$ if it goes to zero
under the restriction map $\H^1(k,A[2])\to \H^1(k_v^{\rm nr},A[2])$,
where $k_v^{\rm nr}$ is the maximal unramified extension of the 
completion $k_v$ of $k$ at $v$. 

\bthe \label{t:main}
Let $k$ be a number field.
Let $A=\prod_{i=1}^r A_i$, where each $A_i$ is a principally
polarised abelian variety. Assume that conditions (a), (b), (c), (d),
(e), (f) and (g) hold. Assume that the $2$-primary subgroup
of the Shafarevich--Tate group $\Sha(A_i^F)\{2\}$ is finite 
for all $i=1,\ldots,r$ and all extensions $F$ of $k$ with 
$[F:k]\leq 2$. Consider the classes in $\H^1(k,A[2])$ that are
unramified at $w_1,\ldots, w_r$ and whose
projection to $\H^1(k,A_i[2])$ is non-zero for each $i=1,\ldots,r$.
If the generalised Kummer variety of $A$ defined by such a class
is everywhere locally soluble, then it has a Zariski dense set 
of $k$-points.
\ethe

\noindent{\bf Remarks} 1. If $r=1$, then condition (d) is not needed 
and condition (e) is vacuous.

\smallskip

\noindent 2. The Brauer--Manin obstruction does not appear 
in the conclusion of the theorem, because the purely algebraic
conditions (a), (b) and (e) imply that the relevant part 
of the Brauer group is trivial, see the remark after Proposition \ref{b7}.

\smallskip

\noindent 
3. If the 2-primary torsion subgroup $\Sha(A_i^F)\{2\}$ is finite, 
then condition (b) implies that the non-degenerate
Cassels--Tate pairing on $\Sha(A_i^F)\{2\}$ is alternating.
See Proposition \ref{f5} based on the work of Poonen--Stoll
\cite{PS} and Poonen--Rains \cite{PR}. 
In the proof of Theorem \ref{t:main} we use a well known consequence of this
result that the number of elements of $\Sha(A_i^F)[2]$ is a square.

\medskip

We employ the following standard notation:

$k_{w_i}$ is the completion of $k$ at $w_i$, 
 
$\O_{w_i}$ is the ring of integers of $k_{w_i}$, 

$\m_{w_i}$ is the maximal ideal of $\O_{w_i}$, and 

$\F_{w_i}=\O_{w_i}/\m_{w_i}$ is the residue field.

\bco \label{f6}
Let $k$ be a number field.
For $i=1,\ldots,r$ let $f_i(x)\in k[x]$ be a monic irreducible 
polynomial of odd degree $n_i\geq 3$ whose Galois group is 
the symmetric group $S_{n_i}$, and let $A_i$ be the Jacobian of 
the hyperelliptic curve $y^2=f_i(x)$.
Assume the existence of distinct odd primes $w_1,\ldots,w_r$ of $k$ such
that $f_i(x)\in\O_{w_j}[x]$ and $\val_{w_i}(\discr(f_j))=\delta_{ij}$
for any $i,j\in\{1,\ldots, r\}$.
Assume that $\Sha(A_i^F)\{2\}$ is finite 
for all extensions $F$ of $k$ with $[F:k]\leq 2$, for all $i=1,\ldots,r$.
Consider the classes in $\H^1(k,A[2])$ that are
unramified at $w_1,\ldots, w_r$ and whose
projection to $\H^1(k,A_i[2])$ is non-zero for each $i=1,\ldots,r$.
If the generalised Kummer variety of $A$ defined by such a class
is everywhere locally soluble, then it has a Zariski dense set 
of $k$-points.
\eco
{\em Proof.} Each $A_i$ is a canonically principally polarised
abelian variety which satisfies
conditions (a), (b), (c), (d) by Lemma \ref{k1}.

Let $C_i$ be the proper, smooth and geometrically integral
curve over $k$ given by the affine
equation $y^2=f_i(x)$, so that $A_i={\rm Jac}(C_i)$. 
As in \cite[Section 4.3]{Liu},
a proper and flat {\em Weierstra{\ss} model} 
$\sC_i$ over $\Spec(\O_{w_i})$ is defined as the normalisation
in $C_i\times_k k_{w_i}$ of the projective line $\P^1_{\O_{w_i}}$ with 
the affine coordinate $x$. Since $2\in\O_{w_i}^*$ the integral closure
of $\O_{w_i}[x]$ in $k_{w_i}(C_i)$ is $\O_{w_i}[x,y]/(y^2-f_i(x))$.
The condition $\val_{w_i}(\discr(f_i))=1$
implies that $\sC_i$ is regular and the special fibre 
$\sC_i\times _{\O_{w_i}} \F_{w_i}$ is geometrically integral with a unique
singular point, which is an ordinary double point, see Cor. 6 and Remark
18 on p. 4602 of \cite{Liu}. In particular, the reduction of $f_i(x)$
modulo $\m_{w_i}$ has one rational double root and $n_i-2$ simple roots.
(This can also be checked directly using Sylvester's formula for the discriminant.)
Now \cite[Thm. 9.6.1]{BLR} implies that the special fibre of the 
N\'eron model of $A_i\times_k k_{w_i}$ over $\Spec(\O_{w_i})$ is connected.
If $j\not=i$, then $\val_{w_i}(\discr(f_j))=0$, and this implies that
$A_j$ has good reduction at $w_i$. We conclude that (f) holds. 

For each $i=1,\ldots,r$ the field $K_i=k(A_i[2])$ is the splitting
field of $f_i(x)$. Since $\Gal(K_i/k)\cong S_{n_i}$, the alternating group
is the unique non-trivial normal subgroup of $\Gal(K_i/k)$. 
Its invariant subfield is $k(\sqrt{\discr(f_i)})$.
Thus if $k'$ is a Galois extension of $k$ 
such that $k\subsetneq k' \subset K_i$, then 
$k(\sqrt{\discr(f_i)})\subset k'$.
The extension $k(\sqrt{\discr(f_i)})$ of $k$ is ramified at $w_i$, 
so (g) holds.

Let $K'_i$ be the compositum of the fields $K_j$ for $j\not=i$.
Since each $K_i$ is a Galois extension of $k$,
the field $K_i\cap K'_i$ is also a Galois extension of $k$.
To verify (e) we need to check that $K_i\cap K'_i=k$ 
for each $i=1,\ldots, r$. Otherwise,
$K_i\cap K'_i$ contains $k(\sqrt{\discr(f_i)})$ 
which is ramified at $w_i$. However, this contradicts
the criterion of N\'eron--Ogg--Shafarevich according to which $K'_i$ 
is unramified at $w_i$, because the abelian varieties $A_j$ for $j\not=i$
have good reduction at this place. Thus (e) holds. $\Box$

\medskip

\noindent{\em Proof of Theorem A.} 
For $i=1,2$ let $C_i$ be the curve of genus 1 given by $y^2=g_i(x)$.
Write $g_i(x)=ax^4+bx^3+cx^2+dx+e$. The classical $\SL(2)$-invariants of
the corresponding quartic binary form $G_i(u,v)=v^4g_i(u/v)$ are
$$I=12ae-3bd+c^2,\quad J=72ace+9bcd-27ad^2-27eb^2-2c^3.$$
Then the Jacobian of $C_i$ is the elliptic curve 
$E_i$ with the equation $u^2=p_i(t)$,
where $p_i(t)=t^3-27Ix-27J$ is the resolvent cubic polynomial of $g_i(x)$,
see \cite[Prop. 3.3.6 (a)]{Sk}. The 0-dimensional scheme
$g_i(x)=0$ is a $k$-torsor $Z_i$ for $E_i[2]$. Then $C_i$ can be viewed as
the twist of $E_i$ by $Z_i$, that is, $C_i=(E_i\times Z_i)/E_i[2]$, 
where $E_i[2]$ acts simultaneously on both factors. 
The antipodal involution acts on $C_i$ by changing the sign of $y$,
so the Kummer surface $\Kum(C_1\times C_2)$
is the minimal desingularisation of the quotient of 
$C_1\times C_2$ by the involution that acts on each component as
$(x,y)\mapsto (x,-y)$. Thus $z^2=g_1(x)g_2(y)$ is indeed an
affine equation of $\Kum(C_1\times C_2)$.

Since the polynomials $g_1(x)$ and $g_2(x)$ have no roots in $k$,
each of the torsors $Z_1$ and $Z_2$ is non-trivial.
The field of
definition $K_i=k(E_i[2])$ of $E_i[2]$ is the splitting field of $p_i(t)$.
Hence the condition $\Gal(g_1)\simeq S_4$ 
implies $\Gal(K/k)=\Gal(p_i)\simeq S_3$, for $i=1,2$. 
Since the disriminant of the quartic $g_i(x)$ is equal to
the discriminant of its resolvent cubic $p_i(t)$,
and $g_i(x)\in\O_{w_j}[x]$ implies $p_i(t)\in\O_{w_j}[t]$,
the primes $w_1$ and $w_2$ satisfy the assumption in Corollary \ref{f6}.
To be in a position to appeal to that corollary we now
show that $Z_i$ is unramified at both $w_1$ and $w_2$. 

Indeed, let $\mathcal Z_{ij}\subset\P^1_{\O_{w_j}}$ be the closed subscheme
given by $G_i(u,v)=0$, where $G_i(u,v)=v^4g_i(u/v)\in\O_{w_j}[u,v]$. 
For $j\not=i$ the discriminant of $G_i(u,v)$ is
a unit in $\O_{w_j}$, thus $\mathcal Z_{ij}$ is 
a finite and \'etale $\O_{w_j}$-scheme
of degree $4$ with the generic fibre $Z_i\times_k k_{w_j}$, 
hence $Z_i$ is unramified at $w_j$. For $i=j$ the discriminant 
of $G_i(u,v)$ is a generator of the maximal ideal of $\O_{w_i}$.
This implies that the fibre $\mathcal Z_{ii}\times_{\O_{w_i}}\F_{w_i}$ 
at the closed point of $\Spec(\O_{w_i})$
is the disjoint union of a double $\F_{w_i}$-point and a 
reduced 2-point $\F_{w_i}$-scheme. The latter gives rise
to two sections of  the morphism
$$\mathcal Z_{ii}\times_{\O_{w_i}}\O_{w_i}^{\rm nr}\to 
\Spec(\O_{w_i}^{\rm nr}).$$
Hence $Z_i$ is unramified at $w_i$.
An application of Corollary \ref{f6} finishes the proof. $\Box$

\medskip

\noindent{\em Proof of Theorem B.} The condition
$\val_w(\discr(f))=1$ implies that $k(\sqrt{\discr(f)})$
has degree 2 over $k$. Hence the Galois group of $f(x)$ 
is not a subgroup of the alternating group $A_5$. 
Any proper subgroup of $S_5$ which acts transitively on $\{1,2,3,4,5\}$
and is not contained in $A_5$, is conjugate to
$\Aff_5=\F_5\rtimes\F_5^*$, the group
of affine transformations of the affine line 
over the finite field $\F_5$, see~\cite{Bur}.
Let us show that this case cannot occur. Indeed,
in the proof of Corollary \ref{f6} we have seen that
the reduction of $f(x)$ modulo $\m_w$ has one rational 
double root and three simple roots, whereas the integral model
defined by $y^2=f(x)$ is regular. It follows that
over the maximal unramified extension of $k_w$
the polynomial $f(x)$ is the product of three linear
and one irreducible quadratic polynomials.
Hence the image of the inertia subgroup in $\Aff_5$
is generated by a cycle of length 2. This is a contradiction
because the elements of order 2 in $\Aff_5$ are
products of two cycles.

We conclude that the Galois group of $f(x)$ is $S_5$. 
The theorem now follows from 
Corollary \ref{f6} provided we check that the relevant class in 
$\H^1(k,A[2])$ is non-zero and unramified at $w$.

For this it is enough to prove that the corresponding $k$-torsor 
for $A[2]$ has no $k$-points but has a $k_w^{\rm un}$-point. 
This torsor is the subset 
$Z\subset {\rm R}_{L/k}(\G_m)/\{\pm 1\}$ given by $\lambda=z^2$.
The natural surjective
map $${\rm R}_{L/k}(\G_m)\to{\rm R}_{L/k}(\G_m)/\{\pm 1\}$$
is a torsor for $\mu_2$. Thus $Z(k)$
is the disjoint union of the images of $k$-points of the torsors
$\lambda=tz^2$ for ${\rm R}_{L/k}(\G_m)$, where $t\in k^*$. Hence
$Z(k)\not=\emptyset$ if and only if $\lambda\in k^*L^{*2}$.
Next, the group $\H^1(k_w^{\rm un},\mu_2)$
consists of the classes of 1-cocycles defined by $1$ and $\pi$,
where $\pi$ is a generator of $\m_w$. Hence $Z(k_w^{\rm un})$
is the disjoint union of the images of $k_w^{\rm un}$-points of the torsors
$\lambda=z^2$ and $\pi\lambda=z^2$ for ${\rm R}_{L/k}(\G_m)$.
We know that for some $\varepsilon\in\{0,1\}$
the valuation of $\pi^\varepsilon\lambda$ 
at each completion of $L$ over $w$ is even. 
Then the torsor for ${\rm R}_{L/k}(\G_m)$ given by $\pi^\varepsilon\lambda=z^2$
has a $k_w^{\rm un}$-point, because any unit is a square as
the residue field
of $k_w^{\rm un}$ is separably closed of characteristic different
from 2. It follows that $Z(k_w^{\rm un})\not=\emptyset$. $\Box$

\section{Galois theory of finite torsors} \label{one}

This section develops some ideas of Mazur and Rubin, see \cite[Lemma 3.5]{MR10}.

We shall work with groups obtained as the semi-direct
product of a group $G$ with a $G$-module $M$, 
where the following lemma will be useful.

\ble \label{f1}
Let $G$ be a group and let $M$ be a simple and faithful $G$-module.
For each normal subgroup $H\subset M^n\rtimes G$ we have
either $M^n\subset H$ or $H\subset M^n$.
\ele
{\em Proof.} Suppose that $M^n$ is not contained in $H$.
The subgroup $H\cap M^n$ is normal in $M^n\rtimes G$,
thus $H\cap M^n$ is a proper $G$-submodule 
of the semisimple $G$-module $M^n$.
Then there is a surjective homomorphism of $G$-modules 
$\varphi:M^n\to M$ such that $\varphi(H\cap M^n)=0$.
The quotient of $M^n\rtimes G$ by $\Ker(\varphi)$ is isomorphic
to $M\rtimes G$. Let us denote by $\rho:M^n\rtimes G\to M\rtimes G$
the corresponding surjective map. The groups $\rho(H)$ and $M$ are normal
in $M\rtimes G$ and $\rho(H)\cap M=\{e\}$, hence they centralise
each other. The image of $\rho(H)$ in $G$ must be trivial,
since $\rho(H)$ acts trivially on $M$, which is a faithful $G$-module.
Thus $\rho(H)=\{e\}$. $\Box$

\medskip

Let $k$ be a field. Let $\bar k$ be a separable closure 
of $k$ and $\Ga_k=\Gal(\bar k/k)$.

Let $M$ be a finite $\Ga_k$-module such that $pM=0$ for a prime 
$p\not={\rm char}(k)$.
Then $M$ is identified with the group of $\bar k$-points
of a finite \'etale commutative group $k$-scheme $\GG_M$.
A cocycle $c:\Ga_k\to M=\GG_M(\bar k)$ gives rise to a twisted action 
of $\Ga_k$ on $\GG_M(\bar k)$, defined as the natural action of 
$\Ga_k$ on $M$ followed by the translation by $c$. The quotient
of $\Spec(\bar k[\GG_M])$ by the twisted action is a $k$-torsor of $\GG_M$.
It comes equipped with a 
$\bar k$-point corresponding to the neutral element of $\GG_M$. 
Conversely, suppose we are given a $k$-torsor $Z$ for $\GG_M$.
For any $z_0\in Z(\bar k)$ the map $c:\Ga_k\to M$
determined by the condition $c(\gamma)^{-1} z_0 = {}^\gamma z_0$ is 
a cocycle $\Ga_k\to M$.
These constructions describe a bijection between $\H^1(k,M)$ and the
set of isomorphisms classes of $k$-torsors for $\GG_M$. 
See \cite[Section 2.1]{Sk}, and also \cite[Ch. 6]{BLR}.

Let $G$ be the finite group $\varphi(\Ga_k)$, where 
$\varphi:\Ga_k\to\Aut\,M$ is the action of $\Ga_k$ on $M$.
Then $M$ is a faithful $G$-module.
Let $K=\bar k^{\Ker(\varphi)}$ so that $\Gal(K/k)=G$.
In other words, $K/k$ is the smallest extension such that 
$\Ga_K$ acts trivially on $M$.
The action $\varphi$ defines a semi-direct product 
$M\rtimes G$. For a cocycle $c:\Ga_k\to M$
the action of $\Ga_k$ on the $\bar k$-points of 
the corresponding torsor $Z$ 
is through the homomorphism $\Ga_k\to M\rtimes G$
given by $(c,\varphi)$. The invariant subfield of the kernel
of this homomorphism is the smallest extension of $k$
such that each $\bar k$-point of $Z$ is defined over this field. 

Let $T$ be a non-empty finite set and let 
$$\alpha:T\to \H^1(k,M)\backslash\{0\}$$
be a map of sets. Let $M^T$ be the direct sum of copies of 
the faithful $G$-module $M$ indexed by $t\in T$. 
For each $t\in T$ let $Z_t$ be a $k$-torsor for the group scheme 
$\GG_M$ such that the class of $Z_t$ in $\H^1(k,M)$ is $\alpha_t$. 
We denote by $K_t\subset\bar k$ the smallest extension
of $k$ such that each $\bar k$-point of $Z_t$ is defined over $K_t$.
Since $Z_t$ is a $k$-torsor for $\GG_M$, each $K_t$ contains $K$. 
Let $K_T$ be the compositum of the fields $K_t$ for $t\in T$.
In other words, $K_T\subset\bar k$ is the smallest Galois extension of $k$
such that $\Gal(\bar k/K_T)$ acts trivially on the $\bar k$-points
of $Z_T=\prod_{t\in T}Z_t$. Write
$$G_T=\Gal(K_T/k)=\Ga_k/\Ga_{K_T}, \quad\quad W_T=\Gal(K_T/K)=\Ga_K/\Ga_{K_T},$$
so that there is an exact sequence
\begin{equation}
1\lra W_T\lra G_T\stackrel{\varphi}\lra G\lra 1.\label{a4}
\end{equation}
For each $t\in T$ the restriction of $\alpha_t$ to $K_T$ is zero,
hence each $\alpha_t$ belongs to the subgroup 
$$\H^1(G_T,M)\subset \H^1(k,M).$$
The class of $Z_T$ as a $k$-torsor for $(\GG_M)^T$ is the direct sum of
the $\alpha_t$, for $t\in T$. We denote this direct sum by $\alpha_T$.
If $a_T:\Ga_k\to M^T$ is a cocycle representing $\alpha_T$,
then $G_T$ acts on $Z_T(\bar k)$ through the 
homomorphism $(a_T,\varphi):G_T\to M^T\rtimes G$. By the definition of $K_T$
this homomorphism is injective.
Since $M$ is a trivial $\Ga_K$-module,
the restriction of $\alpha_T$ to $W_T$ is an injective homomorphism of 
$G$-modules $\tilde\alpha_T:W_T\to M^T$, and we have a commutative diagram
$$\xymatrix{
1\ar[r]& W_T\ar[r]\ar@{^{(}->}[d]_{\tilde\alpha_T}& G_T\ar[r]
\ar@{^{(}->}[d]_{(a_T,\varphi)}& G\ar[r]\ar[d]^{=}& 1\\
1\ar[r]& M^T\ar[r]& M^T\rtimes G\ar[r]& G\ar[r]& 1}$$

When $M$ is a simple $G$-module, the $\F_p$-algebra
$\End_G(M)$ is a division ring by Schur's lemma, hence
is a finite field extension of $\F_p$.

\bpr \label{a1}
Let $M$ be a simple $\Ga_k$-module such that $pM=0$ and $\H^1(G,M)=0$.
Suppose that $\End_G(M)=\F_q$ for some $q=p^s$.
The following conditions are equivalent:

\smallskip

\noindent{\rm (i)}
the classes $\alpha_t$, $t\in T$, are linearly independent in 
the $\F_q$-vector space $\H^1(k,M)$;

\noindent{\rm (ii)}
the classes $\alpha_t$, $t\in T$, form a basis of
the $\F_q$-vector space $\H^1(G_T,M)$;

\noindent{\rm (iii)} $\tilde\alpha_T:W_T\tilde\lra M^T$ 
is an isomorphism of $G$-modules;

\noindent{\rm (iv)} $(a_T,\varphi)$ is an isomorphism of 
groups $G_T\tilde\lra M^T\rtimes G$;

\noindent{\rm (v)} the fields $K_t$, $t\in T$, 
are linearly disjoint over $K$.
\epr
{\em Proof.} Let us prove the equivalence of (i) and (iii).
Consider the inflation-restriction exact sequence
\begin{equation}
0\to \H^1(G,M)\to \H^1(G_T,M)\to \H^1(W_T,M)^G=\Hom_G(W_T,M).
\label{c2}
\end{equation}
The terms of this sequence are $\F_q$-vector spaces via the
$G$-equivariant action of $\F_q$ on $M$.
Since $\H^1(G,M)=0$, the restriction map in (\ref{c2}) is injective. Thus
(i) holds if and only if the homomorphisms $\tilde\alpha_t$, $t\in T$, 
are linearly independent in the $\F_q$-vector space $\Hom_G(W_T,M)$. 

We claim that $\tilde\alpha_T(W_T)$, like any $G$-submodule of $M^T$,
has the form $M\otimes_{\F_q}V$ for some $\F_q$-vector subspace
$V\subset (\F_q)^T$. Indeed, the $G$-module $M$ is simple, 
thus the $G$-module $M^T=M\otimes_{\F_q}(\F_q)^T$ is semisimple,
so any $G$-submodule of $M^T$ is isomorphic
to a direct sum of copies of $M$. Consider a $G$-submodule of $M^T$
isomorphic to $M$. Without loss of generality we can assume that
the projection $\pi_t:M\to M$ is non-zero for each $t\in T$.
Since $M$ is a simple $G$-module, this must be an isomorphism.
For any $t,t'\in T$ the automorphism
$\pi_{t'}\pi_{t}^{-1}:M\to M$ is $G$-equivariant,
so it is the multiplication by a non-zero element of $\F_q$. 
This gives a vector $v\in(\F_q)^T$ such that our submodule 
is $M\otimes_{\F_q}v$. This proves the claim.

The homomorphisms $\tilde\alpha_t$, $t\in T$, 
are linearly independent in the $\F_q$-vector space
$\Hom_G(W_T,M)$ if and only if $V=(\F_q)^T$. Thus (i) and (iii)
are equivalent. 

The equivalence of (iii) and (iv) is clear from the diagram. 

To prove the equivalence of (iii) and (v) recall that
$K_T$ is the compositum of the extensions $K_t/K$, for $t\in T$,
and $W_T=\Gal(K_T/K)$. It is clear that
$\tilde\alpha_T:W_T\to M^T$ is the composed homomorphism 
\begin{equation}
\Gal(K_T/K)\lra\prod_{t\in T}\Gal(K_t/K)\stackrel{\tilde\alpha_t}\lra M^T.
\label{cc}
\end{equation}
Using (\ref{c2}) in the case when $T$ is a one-element set 
one shows that a non-zero class $\alpha_t$ restricts to an isomorphism 
$\tilde\alpha_t:\Gal(K_t/K)\tilde\lra M$.
This shows that (iii) holds if and only if the first arrow in (\ref{cc})
is an isomorphism, which is exactly (v). 

Finally, (ii) obviously implies (i) since $\H^1(G_T,M)\subset \H^1(k,M)$.
On the other hand, we have seen that (i) implies (iii), so we have 
$$\Hom_G(W_T,M)\simeq \Hom_G(M^T,M)=\End_G(M)^T=(\F_q)^T.$$
In view of (\ref{c2}) this shows that the dimension of 
the $\F_q$-vector space $\H^1(G_T,M)$ is at most $|T|$.
Thus if the classes $\alpha_t$, $t\in T$, are linearly independent
in the $\F_q$-vector space $\H^1(G_T,M)$, they form 
a basis of this space. 
We conclude that (ii) is equivalent to (i). $\Box$

\medskip

We record an amusing corollary of this proposition.

\bco
In the assumptions of Proposition \ref{a1} two classes in $\H^1(k,M)$ 
generate the same $\F_q$-vector subpace of $\H^1(k,M)$
if and only if the associated torsors are isomorphic as $k$-schemes.
Moreover, $k$-torsors attached to non-zero classes are connected.
\eco
{\em Proof.} Let $Z_\alpha$ and $Z_\beta$ be $k$-torsors for $\GG_M$
such that the class of $Z_\alpha$ is $\alpha$ and the class of $Z_\beta$
is $\beta$. If $\alpha=0$, then $Z_\alpha$ has a connected component
$\Spec(k)$, so $Z_\alpha\cong Z_\beta$ implies $\beta=0$.
Now assume that $\alpha$ and $\beta$ are non-zero classes in $\H^1(k,M)$
such that $\beta=r\alpha$ for some $r\in\F_q$.
The multiplication by $r$ map $\GG_M\to\GG_M$
descends to an isomorphism of $k$-schemes $Z_\alpha\tilde\lra Z_\beta$
(which is not an isomorphism of torsors unless $r=1$).

Let us prove the converse. Let $K_\alpha$ be the smallest
subfield of $\bar k$ such that each $\bar k$-point of $Z_\alpha$ is defined
over this field, and similarly for $K_\beta$.
If $Z_\alpha\cong Z_\beta$, then clearly $K_\alpha=K_\beta$, in particular,
these field extensions of $K$ are not disjoint. 
By Proposition \ref{a1}, $\alpha$ and $\beta$
generate the same $\F_q$-vector subpace of $\H^1(k,M)$. 

The last statement follows from the equivalence of (ii) and (iii)
in Proposition \ref{a1}, which shows that the base change
from $k$ to $K$ of any $k$-torsor with a non-zero class in $H^1(k,M)$
is a connected $K$-scheme. $\Box$

\medskip

An element $\gamma\in G_T$ defines a map
$$c_\gamma: \H^1(G_T,M)\lra \H^1(\hat\Z,M)=M/(\gamma-1)=M/(g-1)$$
induced by the homomorphism $\hat \Z\to G_T$ that sends $1$
to $\gamma$. Here we denote by $g$ the image of $\gamma$ under the natural
surjective map $G_T\to G$.

\bco \label{a2}
In the assumptions of Proposition \ref{a1} suppose that
$\{\alpha_t|t\in T\}$ is a linearly independent subset
of the $\F_q$-vector space $\H^1(k,M)$.
Fix any $g\in G$ and let $\phi:T\to M/(g-1)$ be a map of sets. 
Then $g$ has a lifting $\gamma\in G_T$
such that $c_\gamma(\alpha_t)=\phi(t)$ for each $t\in T$.
\eco
{\em Proof.} Let $\gamma\in G_T$ be a lifting of $g$.
It is enough to show that the coset $\gamma W_T$ is mapped
surjectively onto $M^T/(g-1)$ by the function that sends
$\gamma x$, $x\in W_T$, to the class of $a_T(\gamma x)$ in
$M^T/(g-1)$. Since $a_T$ is a cocycle we have
$$a_T(\gamma x)=
{}^\gamma a_T(x)+a_T(\gamma)={}^\gamma \tilde\alpha_T(x)+a_T(\gamma)\ 
\in \ M^T.$$
By Proposition \ref{a1} we know that 
$\tilde\alpha_T(W_T)=M^T$, so we are done. $\Box$

\medskip

In the previous statements we did not exclude the case
of the trivial $\Ga_k$-module $M$ (e.g. $M=\F_p$) in which case $G=\{e\}$.
In the next statement we assume that the action of
$\Ga_k$ on $M$ is non-trivial, that is, $G\not=\{e\}$.

\bco \label{f4}
Let $M$ be a simple $\Ga_k$-module such that $pM=0$,
the action of $\Ga_k$ on $M$ is non-trivial and $\H^1(G,M)=0$.
Suppose that $\End_G(M)=\F_q$ for some $q=p^s$, and
$\{\alpha_t|t\in T\}$ is a linearly independent subset
of the $\F_q$-vector space $\H^1(k,M)$.
Then each subfield of $K_T$ which is abelian over $k$ is contained in $K$.
\eco
{\em Proof.} By Proposition \ref{a1} we have $G_T\simeq M^T\rtimes G$.
Let $M_G$ be the module of coinvariants of $M$, and let $N$ be the 
kernel of the canonical surjective map $M\to M_G$. Then $N$
is a $G$-submodule of $M$, which is a subgroup generated by 
${}^\gamma m-m$ for $\gamma\in G$ and $m\in M$. Since $M$ is a simple
$G$-module, we have $N=0$ or $N=M$. But $M$ is a faithful
$G$-module and $G\not=\{e\}$, so the case $N=0$ is not possible.
Hence $N=M$. Note that ${}^\gamma m-m$ is the commutator of $\gamma$
and $m$ in $G_T$. Thus $M^T\subset [G_T,G_T]$, so that
$G_T^{\rm ab}=G^{\rm ab}$, hence the corollary. $\Box$

\medskip

In this paper we shall apply these results to the case when $M$ is the 2-torsion
subgroup $A[2]$ of an abelian variety $A$ over a field $k$,
${\rm char}(k)\not=2$, such that $A$ satisfies conditions (a) and (b)
in Section \ref{2}. Until the end of this section we assume that
$k$ is any field of characteristic different from 2.

\bpr \label{a3}
Suppose that abelian varieties $A_1,\ldots,A_r$ satisfy
conditions (a), (b) and (e). 
Let $Z_i$ be a non-trivial $k$-torsor for $A_i[2]$,
for each $i=1,\ldots, r$, and let $Z=\prod_{i=1}^r Z_i$.
Let $L$ be the \'etale $k$-algebra $k[Z]$, so that $Z\cong\Spec(L)$.
Then $L$ is a field which contains no quadratic extension of $k$.
\epr
{\em Proof.} Let $G_i=\Gal(K_i/k)$, where $K_i$ is the field of 
definition of $A_i[2]$.  The compositum $K=K_1\ldots K_r$
is the field of definition of $A[2]$, where $A=A_1\times\ldots
\times A_r$.
By assumption the fields $K_1,\ldots, K_r$ are linearly
disjoint over $k$, so that the Galois group $G=\Gal(K/k)$ is the product
$G=\prod_{i=1}^r G_i$.

Let $\alpha\in\H^1(k,A[2])$ be the class defining $Z$, so that 
$\alpha=\sum_{i=1}^r\alpha_i$, where each $\alpha_i\in\H^1(k,A_i[2])$ 
is non-zero. For each $i=1,\ldots, r$ we have $A_i[2]^{G_i}=0$
and $\H^1(G_i,A_i[2])=0$ by conditions (a) and (b).
The inflation-restriction exact sequence for $G_i\subset G$ then gives $\H^1(G,A_i[2])=0$.
Therefore we have $\H^1(G,A[2])=0$. 
From (\ref{c2}) with $M=A_i[2]$ 
it follows that each $\alpha_i$ gives rise to a non-zero homomorphism 
$\tilde\alpha_i\in \Hom_G(\Ga_K,A_i[2])$.
The exact sequence (\ref{a4}) takes the form
$$1\to W_\alpha\to G_\alpha\to G\to 1,$$
and we claim that the injective map 
$$\tilde\alpha=\sum_{i=1}^r\tilde\alpha_i:W_\alpha\to A[2]=\oplus_{i=1}^rA_i[2]$$
is an isomorphism of $G$-modules. Indeed, by condition (a) and 
the linear disjointness of
$K_1,\ldots,K_r$ the $G$-modules $A_i[2]$ are simple and pairwise non-isomorphic.
Therefore, if $\tilde\alpha(W_\alpha)\not=A[2]$, then 
$\tilde\alpha_i=0$ for some $i$, which is a contradiction. 

Now we are in a position to complete the proof. 
Let us fix a base point $\bar x\in Z(\bar k)$. Since $\tilde\alpha$
is an isomorphism, each element of $G$ lifts to a unique element of $G_\alpha$
that fixes $\bar x$. This defines a homomorphism $s:G\to G_\alpha$ 
which is
a section of the surjective map $G_\alpha\to G$. Then $s(G)\subset G_\alpha$
is the stabiliser of $\bar x$ in $G_\alpha$. This identifies
$G_\alpha$ with the semi-direct product $A[2]^r\rtimes s(G)$.

It is clear that $G_\alpha$ acts transitively on $Z(\bar k)$ 
since already $W_\alpha$ acts (simply) transitively on $A[2]$. 
It follows that $Z$ is connected,
and hence $L$ is a field. Moreover, $L\cong (K_\alpha)^{s(G)}$,
where $K_\alpha$ is the subfield of $\bar k$ such that
$\Gal(K_\alpha/k)=G_\alpha$.
If $L$ contains a quadratic extension of $k$,
then $s(G)$ is contained in a normal subgroup $H\subset G_\alpha$ 
of index 2. Since $s$ is a section,
the induced homomorphism $H\to G$ is surjective, so its kernel
is a $G$-submodule of $A[2]$ which is a subgroup of $A[2]$ of index 2.
It must be a direct sum of some of the simple
$G$-submodules $A_i[2]$, but then its index is at least $4$. $\Box$

\bco \label{f2}
Suppose that abelian varieties $A_1,\ldots,A_r$ satisfy
conditions (a), (b) and (e). Let $\{\alpha_t|t\in T\}$ 
be a linearly independent subset of
the $\F_2$-vector space $\H^1(k,A_1[2])$. Then the fields
$K_{1,T},K_2,\ldots, K_r$ are linearly disjoint.
\eco
{\em Proof.} Write $E=K_{1,T}\cap K_2\ldots K_r$. In view of condition (e)
it is enough to show that $E=k$. Indeed, $E$ is a Galois
subfield of $K_{1,T}$, so $\Gal(K_{1,T}/E)$ is a normal
subgroup of $G_{1,T}$. By Proposition \ref{a1} we can use Lemma
\ref{f1} with $G=\Gal(K_1/k)$, $M=A_1[2]$ and $n=|T|$.
It follows that $E\subset K_1$ or $K_1\subset E$.
In the first case $E=k$ because $E$ is contained in 
$K_1\cap K_2\ldots K_r=k$, where the equality holds by
condition (e). By the same condition the second case cannot actually occur,
because then $K_1\subset E\subset K_2\ldots K_r$ which 
contradicts the linear disjointness of $K_1,\ldots,K_r$. $\Box$

\section{Kummer map over a local field} \label{local}

Let $A$ be an abelian variety over a local field $k$ of characteristic zero.
The Kummer exact sequence gives rise to a map
$\delta:A(k)\to \H^1(k,A[2])$, called the Kummer map. 
For $x \in A(k)$ choose $\bar x\in A(\bar k)$ such that $2\bar x=x$.
Then $\delta(x)$ is represented by the cocycle that sends
$\gamma\in\Ga_k$ to ${}^\gamma\bar x -\bar x\in A[2]$.

The Weil pairing is a non-degenerate pairing of $\Ga_k$-modules
$A[2]\times A^t[2]\to \Z/2$.
The induced pairing on cohomology followed by the local invariant
of local class field theory gives a non-degenerate pairing 
of finite abelian groups \cite[Cor. I.2.3]{ADT}
$$\H^1(k,A[2]) \times \H^1(k,A^t[2])\lra \Br(k)[2]
\stackrel{\inv}\lra\frac{1}{2}\Z/\Z.$$
The local Tate duality implies that $\delta(A(k))$ and $\delta(A^t(k))$
are the orthogonal complements to each other under this pairing
(see, e.g., the first commutative diagram in the proof of \cite[I.3.2]{ADT}).

When $A$ is principally polarised, we combine the last pairing
with the principal polarisation $A\tilde\lra A^t$
and obtain a non-degenerate pairing
$$\inv(\alpha\cup\beta): \H^1(k,A[2]) \times \H^1(k,A[2])\lra \Br(k)[2]
\stackrel{\inv}\lra\frac{1}{2}\Z/\Z.$$
Then $\delta(A(k))$ is a maximal isotropic subspace of $\H^1(k,A[2])$.
Note that the pairing $\inv(\alpha\cup\beta)$ is also defined for $k=\R$
and the above statements carry over to this case 
\cite[Thm. I.2.13 (a), Remark I.3.7]{ADT}.

\medskip

Let us recall a well known description of $\delta(A(k))$
when $A$ has good reduction.
Let $\kappa$ be the residue field of $k$, and assume
${\rm char}(\kappa)=\ell\not=2$. Then $\delta(A(k))$
is the unramified subgroup 
$$\H^1_{\rm nr}(k,A[2])=\Ker[\,\H^1(\Ga_k,A[2])\to \H^1(I,A[2])\,],$$
where $I\subset\Ga_k$ is the inertia subgroup.
By N\'eron--Ogg--Shafarevich
the inertia acts trivially on $A[2]$,
so that $\H^1_{\rm nr}(k,A[2])=\H^1(\kappa,A[2])$.
The absolute Galois group $\Gal(\bar\kappa/\kappa)=\Ga_k/I$ is 
isomorphic to $\hat \Z$ with the Frobenius element as a topological
generator. Thus we have a canonical isomorphism
\begin{equation}
\delta(A(k))=A[2]/(\Frob-1). \label{loc}
\end{equation}
Since $\hat\Z$ has cohomological dimension 1,
the spectral sequence $$\H^p(\hat \Z,\H^q(I,A[2]))\Rightarrow
\H^{p+q}(k,A[2])$$ gives rise to the exact sequence
$$0\to A[2]/(\Frob-1)\to \H^1(k,A[2])\to \Hom(I,A[2])^\Frob\to 0.$$
The maximal abelian pro-$2$-quotient of $I$ is isomorphic to $\Z_2$,
and $\Frob$ acts on it by multiplication by $\ell$. Thus
$\Hom(I,A[2])=A[2]$ with the natural action of $\Frob$,
so that 
$$\Hom(I,A[2])^\Frob=A[2]^\Frob=\Ker(\Frob-1:A[2]\to A[2]).$$
It follows that the dimension of the $\F_2$-vector space
$A[2]/(\Frob-1)$ equals the dimension of $A[2]^\Frob$,
and therefore
\begin{equation}
\dim\, \H^1(k,A[2])=2\, \dim \, A[2]/(\Frob-1). \label{f3}
\end{equation}

\medskip

If $F/k$ is a quadratic extension, we write $\delta^F:A^F(k)\to \H^1(k,A[2])$
for the Kummer map of $A^F$.
In the rest of this section we summarise some known results relating 
$\delta$, $\delta^F$ and the norm map $\N: A(F)\to A(k)$.

\ble \label{d3}
We have $\delta(\N(A(F)))=\delta(A(k))\cap\delta^F(A^F(k)) \ \subset \ \H^1(k,A[2])$.
\ele
{\em Proof.} Cf. \cite[Prop. 7]{Kramer} or \cite[Prop. 5.2]{MR07}.
Let $\chi: \Ga_k \to \{\pm1\}$ be the quadratic character associated to $F$.
We choose $\sigma \in \Ga_k$ such that $\chi(\sigma) = -1$. 

Suppose that we have $\delta(x)=\delta^F(y)$, where $x\in A(k)$ and 
$y\in A^F(k)$. Using the embedding $A^F(k)\subset A(F)$ we can consider
$y$ as a point in $A(F)$ such that ${}^\sigma y  = -y$. If 
$\bar y\in A(\bar k)$ is such that $2\bar y=y$, then $\delta^F(y)$ is 
represented by the cocycle that sends $\gamma\in\Ga_k$ to 
$$\chi(\gamma)\ {}^{\gamma}\bar y -\bar y={}^\gamma\bar y -
\chi(\gamma)\bar y\in A[2].$$
The equality between classes $\delta(x)=\delta^F(y)$ means that
the corresponding cocycles are equivalent. By changing $\bar y$ if necessary
we can arrange that these cocycles coincide as functions $\Ga_k\to A[2]$.
We deduce that ${}^\gamma(\bar x-\bar y)=\bar x-\chi(\gamma)\bar y$ 
for every $\gamma \in \Gamma_k$.
It follows that $\bar x-\bar y\in A(F)$ and ${}^\sigma(\bar x-\bar y)=
\bar x+\bar y$. Therefore, $x=2\bar x=\N(\bar x-\bar y)$ is a norm from $A(F)$.

Conversely, suppose that $x=\N(z)=z+{}^\sigma z$ for some some $z\in A(F)$.
Set $y={}^{\sigma}z-z$. We claim that $\delta(x)=\delta^F(y)$.
Write $\bar y=\bar x-z$, then $y=2\bar y$. 
We then have $\bar x - \bar y =z$ and $\bar x+\bar y = {}^{\sigma}z$. 
It follows that for every $\gamma \in \Gamma_k$ we have
${}^{\gamma}(\bar x - \bar y) = \bar{x} - \chi(\gamma)\bar{y}$, and hence
$$ {}^{\gamma}\bar x - \bar x = {}^{\gamma}\bar y - \chi(\gamma)\bar y .$$
This implies $\delta(x)=\delta^F(y)$ as desired. $\Box$

\ble  \label{d1}
Let $A$ be a principally polarised abelian variety over $k$ with 
bad reduction such that the number of 
geometric connected components of the N\'eron model
of $A$ is odd. If $F$ is an unramified quadratic extension
of $k$, then $\delta(A(k))=\delta^F(A^F(k))$.
\ele
{\em Proof.} Since $A$ is principally polarised, it is isomorphic
to its dual abelian variety. It follows from
\cite[Prop. 4.2, Prop. 4.3]{Mazur72} that 
the norm map $\N:A(F)\to A(k)$ is surjective.
By Lemma \ref{d3} we see that 
$\delta(A(k))\subset\delta^F(A^F(k))$. Since $F$ is unramified,
the quadratic twist $A^F$ also satisfies the assumptions
of the lemma, and the same argument
applied to $A^F$ gives the opposite inclusion. $\Box$

\ble  \label{d2}
Assume that the residue characteristic of $k$ is not $2$.
If $A$ is an abelian variety over $k$ with good reduction and
$F$ is a ramified quadratic extension
of $k$, then $\delta(A(k))\cap\delta^F(A^F(k))=0$.
\ele
{\em Proof.} In this case we have $\N(A(F))=2A(k)$. 
If $\dim(A)=1$ this is proved in \cite[Lemma 5.5 (ii)]{MR07},
and the same proof works in the general case. 
It remains to apply Lemma \ref{d3}. $\Box$

\section{Selmer group and Cassels--Tate pairing} \label{global}

Let $A$ be an abelian variety over a field $k$ of characteristic zero.
Let $\NS(\ov A)$ be the N\'eron--Severi group of $\ov A$.
The dual abelian variety $A^t$ represents the functor $\Pic^0_A$. 
In particular, we have an exact sequence of $\Ga_k$-modules
\begin{equation}
0\to A^t(\bar k)\to \Pic(\ov A)\to \NS(\ov A)\to 0.
\label{g1}
\end{equation}
The antipodal involution $[-1]:A\to A$  
induces an action $[-1]^*:\Pic(\ov A)\to \Pic(\ov A)$
which turns (\ref{g1}) into
an exact sequence of $\Z/2$-modules. The induced action on 
$\NS(\ov A)$ is trivial. The induced action on $A^t(\bar k)$
is the antipodal involution. Since $A^t(\bar k)$ is divisible,
we obtain $\H^1(\Z/2,A^t(\bar k))=0$.
Thus the long exact sequence of cohomology gives an exact sequence
\begin{equation}
0\to A^t[2]\to \Pic(\ov A)^{[-1]^*}\to \NS(\ov A)\to 0.
\label{g2}
\end{equation}

The group $\NS(\ov A)^{\Ga_k}$ is canonically isomorphic
to the group $\Hom(A,A^t)^{\rm sym}$ of self-dual 
$k$-homomorphisms of abelian
varieties $A\to A^t$. A polarisation on $A$ is an element 
$\lambda\in \NS(\ov A)^{\Ga_k}$. The polarisation is called
principal if the associated morphism $\varphi_\lambda:A\to A^t$
is an isomorphism. Let us write $c_\lambda$ for the image of 
$\lambda$ under the differential $\NS(\ov A)^{\Ga_k}
\to \H^1(k,A^t[2])$ attached to (\ref{g2}). For example,
if $A$ is the Jacobian of a curve $C$ and $\lambda$ is the canonical
principal polarisation of $A$, then $c_\lambda$ is the image
of the class of the theta characteristics torsor of $C$
under the isomorphism $\varphi_\lambda:\H^1(k,A[2])\tilde\lra \H^1(k,A^t[2])$;
see \cite[Thm. 3.9]{PR}.

\ble \label{g3} 
Let $A$ be an abelian variety over $k$ with polarisation $\lambda$.
Then $c_\lambda$ belongs to the kernel of the restriction map 
$\H^1(k,A^t[2])\to \H^1(K,A^t[2])$ for $K=k(A[2])$.
\ele
{\em Proof.} Following Poonen and Rains \cite[Section 2.1]{PR} we associate
to any $\F_2$-vector space $M$ the group $UM$ of invertible elements of
the quotient $\bigwedge M /\bigwedge^{\geq 3}M$, where $\bigwedge M$ is the
exterior algebra of $M$ over $\F_2$.
The homomorphism $UM\to M$ given by the first graded factor gives 
rise to the exact sequence of $\F_2$-vector spaces
\begin{equation}
0\to {\bigwedge}^2 M \to UM\to M\to 0,
\label{g4}
\end{equation}
see \cite{PR}, Remark 2.3 (b).
If $M$ is a $\Ga_k$-module, then (\ref{g4}) is also
an exact sequence of $\Ga_k$-modules.
The exact sequence dual to (\ref{g4}) for $M=A[2]$ fits into the following
commutative diagram:
\begin{equation} \begin{array}{ccccccccc}
0&\to& A^t[2]&\to& \Pic(\ov A)^{[-1]^*}&\to& \NS(\ov A)&\to& 0\\
&&||&&\downarrow&&\downarrow&&\\
0&\to& \Hom(A[2],\Z/2)&\to& \Hom(UA[2],\Z/2)&\to &
\Hom({\bigwedge}^2 A[2],\Z/2)&\to&0
\end{array}
\label{g5}
\end{equation}
see \cite{PR}, Prop. 3.2 and diagram (16). 
The left hand vertical map in (\ref{g5}) is induced by the
Weil pairing $A[2]\times A^t[2]\to \Z/2$. 

Consider (\ref{g5}) as a diagram of $\Ga_K$-modules.
The commutativity of (\ref{g5}) implies that the differential
$\NS(\ov A)^{\Ga_K}\to \H^1(K,A^t[2])$ defined by the upper
row factors through the differential defined by the lower row
$$\Hom({\bigwedge}^2 A[2],\Z/2)^{\Ga_K}\lra \H^1(K,\Hom(A[2],\Z/2)).$$
But this differential is zero, since
$\Ga_K$ acts trivially on $A[2]$ and hence
on all the terms of the lower row of (\ref{g5}). $\Box$

\medskip

Now let $k$ be a number field.
For a place $v$ of $k$ let
$$\loc_v:\H^1(k,A[2])\to \H^1(k_v,A[2])$$
be the natural restriction map. 
If $v$ is a place of good reduction, then
$$\loc_v: \Sel(A)\to A[2]/(\Frob_v-1)$$
is the map provided by (\ref{loc}).

The 2-Selmer group 
$\Sel(A)\subset \H^1(k,A[2])$ is defined as the set
of elements $x$ such that $\loc_v(x)\in\delta(A(k_v))$
for all places $v$ of $k$. Since $A^F[2]=A[2]$ the 2-Selmer groups 
of all quadratic twists $\Sel(A^F)$ are subgroups of $\H^1(k,A[2])$.
We have the well known exact sequence
\begin{equation}
0\to A(k)/2\to \Sel(A)\to \Sha(A)[2]\to 0.\label{lent}
\end{equation}
The Cassels--Tate pairing is a bilinear pairing
$$\langle,\rangle:\Sha(A)\times\Sha(A^t)\to\Q/\Z.$$
If $\Sha(A)$ is finite, then $\Sha(A^t)$ is finite too and the Cassels--Tate
pairing is non-degenerate, see \cite{ADT}. A polarisation $\lambda$
on $A$ induces a homomorphism ${\varphi_\lambda}_*:\Sha(A)\to \Sha(A^t)$.

\bpr \label{f5}
Let $A$ be an abelian variety over a number field $k$ with
a principal polarisation $\lambda$. Then
condition (b) implies that the Cassels--Tate pairing 
$\langle x,{\varphi_\lambda}_* y\rangle$ on $\Sha(A)\{2\}$ is alternating.
In particular, if the $2$-primary subgroup $\Sha(A)\{2\}$ is finite, 
then the cardinality of $\Sha(A)[2]$ is a square.
\epr
{\em Proof.} By a result of Poonen and Stoll we know that
$c_\lambda\in \Sel(A^t)$, see \cite[Cor. 2]{PS}.
If $c'_\lambda$ is the image of $c_\lambda$ in
$\Sha(A^t)[2]$, then \cite[Thm. 5]{PS} says that
$\langle x,{\varphi_\lambda}_* x+c'_\lambda\rangle=0$
for any $x\in \Sha(A)$. Thus it is enough to
prove that $c_\lambda=0$. Lemma \ref{g3} implies that
$c_\lambda$ belongs to the image of the inflation map
$\H^1(G,A^t[2])\to \H^1(k,A^t[2])$, where $G=\Gal(k(E[2])/k)$
is the image of $\Ga_k\to\GL(A[2])$. Since
$\lambda$ is a principal polarisation, $\varphi_\lambda$
induces an isomorphism of $\Ga_k$-modules 
$A[2]\tilde\lra A^t[2]$. Now condition (b) implies
$\H^1(G,A^t[2])=\H^1(G,A[2])=0$, hence $c_\lambda=0$. $\Box$

\section{Generalised Kummer varieties} \label{fibr}

Let $A$ be an abelian variety over a field $k$ of characteristic
different from 2. Let $Z$ be a $k$-torsor for the group scheme $A[2]$.
The associated a 2-covering $f:Y\to A$ is a $k$-torsor
for $A$ defined as the quotient of $A\times_k Z$ by the diagonal action of
$A[2]$. In other words, $Y$ is a twisted form of $A$
with respect to the action of $A[2]$ by translations.
The morphism $f$ is induced by the first projection, and we have
$Z=f^{-1}(0)$.

Let $\tilde Y$ be the blowing-up of $Z$ in $Y$. The antipodal involution
acts on $Y$ fixing $Z$ point-wise, and extends to an involution
$\iota:\tilde Y\to \tilde Y$ whose fixed point set is precisely
the exceptional divisor. It is easy to see that the quotient 
$\Kum(Y)=\tilde Y/\iota$ is smooth. 
It is called the {\em generalised Kummer variety} attached
to $A$ and $Z$. We note that the branch locus of $\tilde Y\to \Kum(Y)$ is
$Z\times_k\P^{d-1}_k$, where $d=\dim(A)$. 

Let $F$ be an extension of $k$ of degree at most 2.
The antipodal involution of $A$ 
commutes with translations by the elements of $A[2]$, so
the quadratic twist $Y^F$ of $Y$ is a $k$-torsor for $A^F$.
We have a natural embedding $i_F:Z\to Y^F$.
Then $\tilde Y^F$, defined as the blowing-up of $i_F(Z)$ in $Y^F$, 
is the quadratic twist of $\tilde Y$ that preserves
the projection to $\Kum(Y)$. It is clear that $Y^F$, and hence
$\Kum(Y)$, has a $K$-point over any extension of $k$ such that $\alpha$
is in the kernel of the natural map $\H^1(k,A[2])\to\H^1(K,A^F)$.

\bpr \label{b7}
Let $k$ be a number field, and let
$A=\prod_{i=1}^r A_i$ be a product of abelian varieties over $k$ satisfying
conditions (a), (b), (e) and (f). Let $Z$ be a $k$-torsor for $A[2]$
unramified at $w_1,\ldots,w_r$,
and let $Y$ be the attached $2$-covering of $A$.
If $\Kum(Y)$ is everywhere locally soluble, then
there exists an extension $F$ of $k$ of degree at most $2$
such that $Y^F$ is everywhere locally soluble
and $F$ is split at $w_1,\ldots,w_r$.
\epr
{\em Proof.} Conditions (a), (b) and (e) are satisfied,
so we can use Proposition \ref{a3}. It says that the \'etale $k$-algebra
$L=k[Z]$ is a field and the kernel of the natural 
homomorphism $k^*/k^{*2}\to L^*/L^{*2}$ is trivial.
In this situation \cite[Thm. 3]{SSD} together with
\cite[Lemma 6]{SSD} imply the existence of $F$ such that $Y^F$
is everywhere locally soluble. More precisely,
\cite[Lemma 6 (b)]{SSD} shows that we can arrange for $F$ to be split at
a given finite set of places of $k$, at each of which 
$Y$ has local points. Thus to finish the proof we need to show that
$Y$ has $k_w$-points, where $w$ is any of the places $w_1,\ldots,w_r$.

By assumption the class $[Z]\in\H^1(k,A[2])$ 
goes to zero under the composed map 
$$\H^1(k,A[2])\lra \H^1(k_w,A[2])\lra\H^1(k_w^{\rm nr},A[2]).$$
Hence the class $[Y]\in\H^1(k,A)[2]$ goes to zero under the composed map
\begin{equation}
\H^1(k,A)\lra \H^1(k_w,A)\lra\H^1(k_w^{\rm nr},A). \label{non-r}
\end{equation} 
The second arrow in (\ref{non-r}) is the restriction map
$\H^1(\Ga_{k_w},A)\to \H^1(I_w,A)$, where $\Ga_{k_w}=\Gal(\ov{k_w}/k_w)$
and $I_w\subset \Ga_{k_w}$ is the inertia subgroup.
By the inflation-restriction sequence we see that the class
$[Y\times_k k_w]\in \H^1(\Ga_{k_w},A)$ belongs to the subgroup
$\H^1(\Ga_{k_w}/I_w,A(k_w^{\rm nr}))$. 
Let $\sA\to \Spec(\O_w)$ be the N\'eron model of $A\times_k k_w$.
By \cite[Prop. I.3.8]{ADT} we have an isomorphism
$$\H^1(\Ga_{k_w}/I_w,A(k_w^{\rm nr}))=\H^1(\Ga_{k_w}/I_w, \pi_0(\sA\times_{\O_w}\F_w)),$$
where $\pi_0(\sA\times_{\O_w}\F_w)$
is the group of connected components of the special fibre
of $\sA\to \Spec(\O_w)$.
Since $2[Y]=0$, condition (f) implies that $[Y\times_k k_w]=0$,
hence $Y$ has a $k_w$-point. $\Box$

\medskip

\noindent{\bf Remark} The proof of the existence of $F$ such that $Y^F$
is everywhere locally soluble in \cite[Lemma 6 (b)]{SSD} is an
application of (the proof of)
\cite[Thm. A]{CS}. The latter result allows one
to find an everywhere locally soluble fibre of a projective morphism $\pi:X\to\P^1_k$
which is smooth over $\G_{m,k}\subset\P^1_k$ and has 
the geometrically integral generic fibre, whenever $X$ is
everywhere locally soluble and 
the vertical Brauer group of $\pi:X\to\P^1_k$ creates no obstruction.
In the situation of \cite[Lemma 6]{SSD} the variety $X$
is a smooth compactification of the total space of the
family of all quadratic twists of $Y$ parameterised by $\G_{m,k}$.
By \cite[Thm. 3]{SSD} the injectivity of 
the natural map $k^*/k^{*2}\to L^*/L^{*2}$ implies that
the vertical Brauer group of $\pi:X\to\P^1_k$ comes from $\Br(k)$ 
and hence produces no obstruction.

\section{Proof of Theorem \ref{t:main}} \label{proof}

Suppose that our generalised Kummer variety is $\Kum(Y)$,
where $Y$ is the $k$-torsor for $A$ defined by a class
$\alpha\in \H^1(k,A[2])$. To prove the existence of a $k$-point in $\Kum(Y)$
it is enough to find a quadratic extension $F$
of $k$ such that $\alpha$ goes to 0 in $\H^1(k,A^F)$.
We write $\alpha=\sum_{i=1}^r\alpha_i$, where $\alpha_i\in\H^1(k,A_i[2])$
is non-zero for each $i=1,\ldots, r$.

By Proposition \ref{b7} there is a quadratic extension $F$ of $k$
split at $w_1,\ldots,w_r$ such that $\alpha\in\Sel(A^F)$. 
Replacing $A$ with $A^F$ we can assume without loss of generality
that $\alpha\in \Sel(A)$. By doing so we preserve conditions
(a), (b), (c), (d), (e) and (g) that are not affected by quadratic twisting.
The extension $F/k$ is split at $w_1,\ldots,w_r$, so replacing $A$ 
by $A^F$ also preserves condition (f).

Let $S_0$ be the set of places of $k$ that contains all
the Archimedean places and the places above $2$.

\ble \label{e1}
Let $S$ be the set of places of $k$ which is the union of $S_0$ and 
all the places of bad reduction of $A$ excluding $w_1,\ldots,w_r$.
For each $i=1,\ldots, r$ let $\alpha_i\in \Sel(A_i)$
be non-zero. Let $\beta\in\Sel(A_1)$
be a non-zero class such that $\beta\not=\alpha_1$.
Then there exists $q\in k^*$ such that $\q=(q)$ is a prime ideal of $k$ 
with the following properties:

\begin{enumerate}
\item
all the places in $S$ (including the Archimedean places) 
are split in $F=k(\sqrt{q})$, in particular, $\q\notin S$;
\item
$A$ has good reduction at $\q$;
\item
$\Frob_\q$ acts on $A_1[2]$ as $g$;
\item
$\Frob_\q$ acts on $A_i[2]$ as $h$, for each $i\not=1$;
\item
$\loc^1_\q(\alpha_1)=0$, but $\loc^1_\q(\beta)\not=0$.
\end{enumerate}
\noindent Here the elements $g$ and $h$ are as in conditions (c) and (d).
\ele
{\em Proof.} We adapt the arguments from
the proof of \cite[Prop. 5.1]{MR10}.
Let $T=\{\alpha_1,\beta\}$ and let $K_{1,T}$ be the field
defined in Section \ref{one}.
By Corollary \ref{a2} applied to $M=A_1[2]$ we can find a lift 
$g_T\in G_T=\Gal(K_{1,T}/k)$ of 
$g\in G_1$ such that the associated map
$$c_{g_T}:\H^1(k,A_1[2])\lra A_1[2]/(g-1)\simeq\Z/2$$
annihilates $\alpha_1$ but not $\beta$.

The fields $K_{1,T},K_2,\ldots,K_r$ are Galois extensions
of $k$ that are linearly disjoint by
condition (e) and Corollary \ref{f2}. Let $\K$ be the compositum
of $K_{1,T},K_2,\ldots,K_r$. This is a Galois extension of $k$
with the Galois group $\Gal(\K/k)=G_T\times\prod_{i=2}^r G_i$.

Let the modulus $\m$ be the formal product of the real places of $k$,
$8$ and all the odd primes in $S$. Let $k_\m$ be the ray class field
associated to the modulus $\m$. This is an abelian extension
of $k$ which is unramified away from $\m$. We claim that
$k_\m$ and $\K$ are linearly disjoint over $k$. Indeed, $k'=k_\m\cap\K$
is a subfield of $\K$ that is abelian over $k$ and unramified at
$w_1,\ldots,w_r$. We note that 
$\Gal(\K/k)^{\rm ab}=G_T^{\rm ab}\times\prod_{i=2}^r G_i^{\rm ab}$.
By Corollary \ref{f4} applied to $M=A_1[2]$
we have $G_T^{\rm ab}=G_1^{\rm ab}$.
Therefore, $\Gal(\K/k)^{\rm ab}=\prod_{i=1}^r G_i^{\rm ab}$, so that
$k'$ is contained in the compositum $L=k_1^{\rm ab}\ldots k_r^{\rm ab}$ 
of linearly disjoint abelian extensions
$k_1^{\rm ab},\ldots, k_r^{\rm ab}$. 

Write $M=k_1^{\rm ab}\ldots k_{r-1}^{\rm ab}$.
The extension $k_r^{\rm ab}/k$ is totally ramified at $w_r$
by condition (g), whereas $k'/k$ and $M/k$ are unramified at $w_r$
(the latter by the criterion of N\'eron--Ogg--Shafarevich). Hence
$L/M$ is totally ramified at each prime $v$ of $M$
over $w_r$. 
Since $M\subset k'M\subset L$, where $k'M/M$ is unramified over $v$,
we must have $k'\subset M$. Continuing by induction we prove that
$k'=k$, as required.

It follows that $k_\m\K$ is a Galois extension of $k$
with the Galois group 
$$\Gal(k_\m\K/k)=\Gal(k_\m/k)\times G_T\times\prod_{i=2}^r G_i.$$
By Chebotarev density theorem we can find a place $\q$
such that the corresponding Frobenius element in $\Gal(k_\m\K/k)$
is the conjugacy class of $(1,g_T,h,\ldots,h)$. Then $\q$ is a principal
prime ideal with a totally positive generator $q\equiv 1\bmod 8$,
hence $q$ is a square in each completion of $k$ at a prime over $2$.
We also have $q\equiv 1\bmod \mathfrak p$ for any odd $\mathfrak p\in S$.
Thus all the places of $S$ including the Archimedean places
are split in $F=k(\sqrt{q})$. All other conditions are satisfied 
by construction. $\Box$

\bpr \label{p:core}
For any $\beta\in \Sel(A_1)$, $\beta\not=0$, $\beta\not=\alpha_1$, 
there exists a quadratic extension $F/k$ unramified at
the places of $S_0$ and all the places of bad reduction of $A$,
such that 
$$\Sel(A_1^F) \subset \Sel(A_1), \ 
\alpha_1\in\Sel(A_1^F),\ \beta\notin\Sel(A_1^F),
\ \Sel(A_i^F)=\Sel(A_i) \ \text{for}\ i\not=1.$$
\epr
{\em Proof.} Let $F=k(\sqrt{q})$ be as in Lemma \ref{e1}.
Let $i\in\{1,\ldots, r\}$.
Since $F$ is split at each $v\in S$ we have 
$A_i^F\times_k k_v\cong A_i\times_k k_v$, so that
the Selmer conditions at $S$ are identical for $A_i$ and $A_i^F$.
These conditions are also identical for all primes
where both $A_i$ and $A_i^F$ have good reduction,
and this includes the primes $w_j$ if $j\not=i$.
At $w_i$ the extension $F/k$ is unramified, and
by condition (f) we can apply Lemma \ref{d1}, so we obtain
$\delta(A_i(k))=\delta^F(A_i^F(k))$.

It remains to check the behaviour at $\q$,
which is a prime of good reduction for $A_i$.
If $i\not=1$ then $\Frob_\q=h$, and from condition (d) and formula (\ref{f3})
we deduce $\H^1(k_\q,A_i[2])=0$ so the Selmer conditions for both
$A_i$ and $A_i^F$ at $\q$ are vacuous. This proves that
$\Sel(A_i^F)=\Sel(A_i)$ whenever $i\not=1$.

In the rest of the proof we work with $A_1$. 
The Selmer conditions
for $A_1$ and $A_1^F$ are the same at each place $v\not=\q$.
Thus $\loc^1_\q(\alpha_1)=0$ implies $\alpha_1\in \Sel(A_1^F)$.
Moreover, $\delta(A_1(k_\q))\cap\delta^F(A_1^F(k_\q))=0$
by Lemma \ref{d2}, hence $\beta\not\in\Sel(A_1^F)$.

To prove that $\Sel(A_1^F)\subset\Sel(A_1)$ it is enough to show that
for $A_1$ the Selmer condition at $\q$ is implied by the Selmer conditions
at the other places of $k$. Indeed, let $\xi\in\H^1(k,A_1[2])$
be an element satisfying the Selmer condition at each place $v\not=\q$,
but not necessarily at $\q$.
By global reciprocity the sum of $\inv_v(\beta\cup\xi)\in\frac{1}{2}\Z/\Z$
over all places of $k$, including the Archimedean places, 
is $0$. Since the images of $\xi$ and $\beta$ in $\H^1(k_v,A_1[2])$
belong to $\delta(A_1(k_v))$ for all $v\not=\q$ we obtain 
$\inv_v(\beta\cup\xi)=0$. By the global reciprocity we deduce
$\inv_\q(\beta\cup\xi)=0$. The non-zero element $\loc^1_\q(\beta)$
generates $\delta(A_1(k_\q))$, because
$$\delta(A_1(k_\q))=A_1/(\Frob_q-1)=A_1/(g-1)=\Z/2,$$
where we used (\ref{loc}) and the fact
that $\Frob_\q$ acts on $A_1[2]$ as the element $g$
of condition (c). Since $A_1$ is principally polarised, 
$\delta(A_1(k_\q))$ is a maximal isotropic subspace of 
$\H^1(k_\q,A_1[2])$, therefore $\inv_\q(\beta\cup\xi)=0$
implies that the image of $\xi$ in $\H^1(k_\q,A_1[2])$
lies in $\delta(A_1(k_\q))$. $\Box$

\medskip

\noindent{\em End of proof of Theorem}~\ref{t:main}.
The extension $F/k$ is unramified at all the places where $A$ has
bad reduction, so replacing $A$ by $A^F$ preserves condition (f).
Conditions (a), (b), (c), (d), (e) and (g) 
are not affected by quadratic twisting.
By repeated applications of Proposition~\ref{p:core} we can find 
a quadratic extension $F/k$
such that $\alpha_i$ is the only non-zero element 
in $\Sel(A_i^F)$, for all $i=1,\ldots, r$. 
The exact sequence (\ref{lent}) for $A_i^F$ shows that
$\Sha(A_i^F)[2]$ is of size at most $2$.
If the 2-primary subgroup of $\Sha(A_i^F)$ is finite, then,
by Proposition \ref{f5},
the number of elements in $\Sha(A_i^F)[2]$ is a square.
Thus $\Sha(A_i^F)[2]=0$, so that the image of $\alpha_i$ in $\H^1(k,A_i^F)$
is $0$. Then the image of $\alpha$ in $\H^1(k,A^F)$
is $0$, so that $Y^F\cong A^F$ and hence 
$Y^F(k)\not=\emptyset$. This implies that $\Kum(Y)=\Kum(Y^F)$
has a $k$-point. 

It remains to prove that $k$-points are Zariski dense in $\Kum(Y)$.
For each $i$ the exact sequence 
(\ref{lent}) for $A_i^F$ shows that $A_i^F(k)/2\not=0$.
Since $A_i^F[2](k)=0$ by condition (a), we see that $A_i^F(k)$
is infinite. The neutral connected component of the 
Zariski closure of $A_i^F(k)$ in $A_i^F$ is an abelian subvariety 
$B\subset A_i^F$ of positive dimension. By condition (a) we must
have $B=A_i^F$. Thus the set $A_i^F(k)$ is Zariski dense in $A_i^F$ for each
$i=1,\ldots, r$, so that $A^F(k)$ is Zariski dense in $A^F$. It follows that
$k$-points are Zariski dense in $\Kum(Y)=\Kum(Y^F)\cong\Kum(A^F)$. $\Box$

\noindent 
Department of Mathematics, South Kensington Campus, Imperial College, 
London, SW7 2BZ England, United Kingdom, and

\smallskip

\noindent
Institute for the Information Transmission Problems, 
Russian Academy of Sciences, 19 Bolshoi Karetnyi, Moscow, 127994 Russia

\smallskip

\noindent{\tt a.skorobogatov@imperial.ac.uk}

\smallskip

\noindent 
D\'epartement de Math\'ematiques et Applications,
\'Ecole Normale Sup\'erieure,
45 rue d'Ulm,
Paris 75005,
France

\smallskip

\noindent{\tt harpaz@dma.ens.fr}

\end{document}